\documentclass[conference]{IEEEtran}
\usepackage{cite}

\ifCLASSINFOpdf
   \usepackage[pdftex]{graphicx}
\usepackage{caption}
\usepackage{subcaption}
\else
\fi
\usepackage[cmex10]{amsmath}
\usepackage{amssymb}
\begin{document}
%
\title{Synchronization of Heterogeneous Kuramoto Oscillators with Graphs of Diameter Two}

\author{\IEEEauthorblockN{Andrey Gushchin}
\IEEEauthorblockA{Center for Applied Mathematics\\
Cornell University\\
Ithaca, New York 14850\\
Email: avg36@cornell.edu}
\and
\IEEEauthorblockN{Enrique Mallada}
\IEEEauthorblockA{ Center for the Mathematics \\of Information\\
California Institute of Technology\\
Pasadena, California 91125\\
Email: mallada@caltech.edu}
\and
\IEEEauthorblockN{Ao Tang}
\IEEEauthorblockA{ School of Electrical and \\Computer Engineering\\
Cornell University\\
Ithaca, New York, 14850\\
Email: atang@ece.cornell.edu}}


%


\maketitle

\begin{abstract}
In this article we study synchronization of Kuramoto oscillators with heterogeneous frequencies, and where underlying topology is a graph of diameter two. When the coupling strengths between every two connected oscillators are the same, we find an analytic condition that guarantees an existence of a Positively Invariant Set (PIS) and demonstrate that existence of a PIS suffices for frequency synchronization. For graphs of diameter two, this synchronization condition  is significantly better than existing general conditions for an arbitrary topology.
 If the coupling strengths can be different for different pairs of connected oscillators, we formulate an optimization problem that finds sufficient for synchronization coupling strengths such that their sum is minimal.
\end{abstract}


%
\IEEEpeerreviewmaketitle


\section{Introduction}

Synchronization of coupled oscillators is an important topic of research for scientists from different areas including neuroscience~\cite{c10,c9,c11}, physics~\cite{c13,c12},  mathematics~\cite{c5}, and engineering \cite{c24}, \cite{c26}, \cite{c27}. Kuramoto model \cite{c6} of coupled oscillators, despite its seeming simplicity, demonstrates a quite rich dynamic behavior and has become a canonical model for studying synchronization.

The two main features that describe the behavior of a system of coupled oscillators are the coupling function and the interconnection topology. 
In the case of the Kuramoto model a trigonometric $\sin()$ is used as the coupling function; a broader class of the coupling functions, however,  has also been discussed~\cite{c23,c7,c22, plastic}. The most popular assumption on the interconnection topology is that all oscillators are connected to each other, which corresponds to a fully connected graph or a graph of diameter one~\cite{c1,c2}.
A much more general approach is to study the systems of oscillators with an arbitrary underlying topology~\cite{c20,c3, c28, c5, c19}.

Several additional assumptions can be made to make analysis of the Kuramoto model more tractable. First, one may consider a limit case when the model contains infinite number of oscillators \cite{c4,c6,c14,c8}.
Second, it can be assumed that all oscillators have equal intrinsic frequencies, and therefore form a gradient system of homogeneous oscillators~\cite{c7}. Alternatively, as we do in this article, one may let the frequencies to take distinct values and thus analyze a system of heterogeneous oscillators~\cite{c15,c1,c2,c28,c18,c17,c16}.
Finally, the coupling strengths can be equal for all pairs of connected oscillators, or are allowed to take different values for different connections.

In this paper, we consider a system of finite number of heterogeneous Kuramoto oscillators in which the underlying topology is a graph of diameter two, a natural step to further generalization of the complete graph (diameter one) case. First, we consider the case when the coupling strength is the same for all pairs of connected oscillators and formulate an analytic condition that guarantees boundedness of the trajectories, which in our case also implies synchronization.
While there exist more general synchronization conditions \cite{c20,c3, c28, c5, c19} that are applicable to the systems with an arbitrary topology, they are significantly more restrictive (when applied to the diameter two graphs) compared to our analytic condition. We provide simulation results that illustrate the improvement over existing results for the graphs of diameter two.
Second, when the coupling strengths are allowed to be different for different pairs of interconnected oscillators, we formulate an optimization problem that finds the coupling strengths such that their sum is minimal while synchronization is preserved.

The rest of the paper is organized as follows: in Section II we describe the problem setup  as well as the main challenge for guaranteeing synchronization, i.e. showing boundedness of trajectories. This challenge is addressed in Section III-A, where we show a general, yet hard to check, condition for synchronization (Proposition 1). This condition is made tractable in Section III-B for the case of equal coupling strengths.
Further, in Section III-C we present an optimization approach to study the case when the coupling strengths can be different for different pairs of connected oscillators.
We illustrate our findings using simulations in Section IV and conclude in Section V.

\section{Problem Formulation}

In this article we study a system of Kuramoto oscillators in which each oscillator is described by the following equation:
\begin{equation}
\label{eq:sys1}
\dot \phi_{i}=\omega_{i}+\sum_{j \in N_{i}} \frac{K_{ij}}{n}\cdot \sin(\phi_{j}-\phi_{i}),
\end{equation}
where $N_{i}$ is a set of oscillators connected to oscillator $i$, i.e. the set of its neighbors, $K_{ij}$ is the coupling strength between oscillators $i$ and $j$, and $n$ is the total number of oscillators in the system.
The coupling strength is symmetric ($K_{ij}=K_{ji}, \; \; \forall i,j$), and can be the same for all connections as assumed in Section III-B, or can be different for different pairs of connected oscillators as in Section III-C. 
 We also assume that the intrinsic frequencies of oscillators $\omega_{i}$ are heterogeneous, which implies that they are not necessary equal. Frequencies, however, do not change their values with time, so each $w_{i}$ is a constant.

 In this article we study frequency synchronization of the system \eqref{eq:sys1}. System \eqref{eq:sys1} achieves synchronization if $\dot \phi_{1}(t) = \dot \phi_{2}(t) = \dots = \dot \phi_{n}(t)$ as $t \rightarrow \infty $. We will denote the common phase velocity by $\dot \phi$.

This common phase velocity $\dot \phi$ is an average intrinsic frequency of the oscillators:
\begin{equation*}
\dot \phi = \frac{\sum\limits_{k=1}^{n}\omega_{i} }{n}.
\end{equation*}
Indeed, when $\dot \phi_{1} = \dot \phi_{2} = \dots = \dot \phi_{n}$,  the sum of all the equations of \eqref{eq:sys1} is: $(\dot \phi_{1} + \dot \phi_{2} +\dots +\dot \phi_{n} ) = \omega_{1}+ \omega_{2}+ \dots+~\omega_{n}$, because each $\frac{K_{ij}}{n} \sin(\phi_{j}-\phi_{i})$ is added to $\frac{K_{ji}}{n}\sin(\phi_{i}-\phi_{j})$ and gives zero. 


We now define the deviations of the frequencies for all $i=1,\dots, n$ by $\bar{\omega}_{i}\triangleq \omega_{i}-\bar{\omega}$, where $\bar{\omega}\triangleq~\frac{1}{n}\sum\limits_{k=1}^{n}\omega_{i}$ is the average intrinsic frequency, and study the following system:
\begin{equation}
\label{eq:sys2}
\dot \phi_{i}=\bar{\omega}_{i}+\sum_{j \in N_{i}}\frac{K_{ij}}{n} \cdot \sin(\phi_{j}-\phi_{i}). 
\end{equation}

Each limit cycle of system \eqref{eq:sys1} is an equilibrium of \eqref{eq:sys2}. Therefore, we will focus on finding conditions when system \eqref{eq:sys2} synchronizes, i.e. when $\dot{\phi_{i}}=0 \; \forall i=1,\dots,n$. Due to the rotational invariance of system \eqref{eq:sys2}, and since the phase average remains the same  ($\dot{\phi_{1}}+\dots+\dot{\phi_{n}}=0$), we can assume without loss of generality that $\sum\limits_{i=1}^{n}\phi_{i}^{t}=0$ for all $t\geq0$, where $\vec{\phi}^{t}$ are the trajectories of system \eqref{eq:sys2}.

We will show frequency synchronization of system \eqref{eq:sys2} by providing a Lyapunov function and using LaSalle's Invariance Theorem \cite{c21}.
 When the oscillators are homogeneous, all the intrinsic frequencies are equal, i.e. deviations $\bar{\omega}_{1}= \dots= \bar{\omega}_{n}=0$, and the following Lyapunov function can be used:
\begin{equation*}
V_{0}(\vec \phi)= -\sum_{ij \in E, i<j}\frac{K_{ij}}{n} \cdot \cos(\phi_{i}-\phi_{j}),
\end{equation*}
where $\vec{\phi} \in \mathbb R^{n}$ and $E$ is the edge set of a given graph.
It can be verified that:
\begin{equation*}
\dot{V_{0}}(\vec \phi)=-\sum\limits_{i=1}^{n}\dot{\phi}_{i}^{2}\leq 0
\end{equation*}
Since function $V_{0}(\vec \phi)$ is well-defined on a n-dimensional torus $\mathbb T^{n}$ which is compact, applying the LaSalle's Invariance Theorem (on $\mathbb T^{n}$) guarantees synchronization of the oscillators.

When the intrinsic frequencies are not equal, we have a system of heterogeneous oscillators, and we still can provide a potential function for this case:
\begin{equation}
\label{eq:V}
V(\vec \phi)\triangleq -\sum_{k=1}^{n} (\bar{\omega}_{k}\phi_{k})-\sum_{ij \in E, i<j} \frac{K_{ij}}{n}\cdot(\cos(\phi_{i}-\phi_{j})).
\end{equation}
We can check that the derivative of this function is also non-positive and is equal to zero only at equilibrium, i.e. when the frequencies are synchronized.

The key problem here is that function $V(\vec \phi)$ is not bounded from below and cannot be defined on $\mathbb T^{n}$. Therefore, we are not able to apply directly the LaSalle's Invariance Theorem. However, if we show that the trajectories $\vec{\phi} \in \mathbb R^{n}$ of \eqref{eq:sys2} are bounded, then the function $V(\vec \phi)$ is bounded as well, hence synchronization follows.

One of the techniques for showing boundedness of the trajectories is to find a bounded Positively Invariant Set (PIS) for the oscillators' phases. 
The goal of this article is to show that when some conditions are met, a PIS exists, and if the initial phases are in this PIS, then the trajectories will be bounded and, therefore, system \eqref{eq:sys2} will achieve frequency synchronization.


\section{Main Results}
\label{sec-results}
This section is organized as follows: we first introduce the notations used in this article and provide a general synchronization condition in Proposition 1. We also demonstrate by means of an example that existence of an equilibrium does not guarantee that system \eqref{eq:sys2} achieves frequency synchronization for all initial phase values. In Subsection $B$ we provide an analytic synchronization condition for system \eqref{eq:sys2} with equal coupling strengths. In Subsection $C$ we study a more general case when the coupling strengths can be different for different edges.

\subsection{Preliminary Results}

Let $G=(V, E)$ be an undirected graph with vertex set $V$ and edge set $E$ that defines the topology of the system \eqref{eq:sys2}. Distance between vertices $i$ and $j$ is defined as a number of edges in the shortest path between $i$ and $j$, where the length of a path is defines as the number of edges in it. Diameter of a graph is defined as the maximum distance between its two vertices. All the results presented in this article are formulated for the graphs of diameter two.

We denote by $A$ the symmetric adjacency matrix of a graph $G$, and define for each pair of vertices $i, j$ constant $P_{ij}$:
\begin{equation}
P_{ij}\triangleq {\bf{a_{i}}}\cdot{\bf{a_{j}^{T}}} + 2A_{ij},
\end{equation}
where ${\bf{a_{i}}}$ and ${\bf{a_{j}}}$ are the $i^{th}$ and $j^{th}$ rows of matrix $A$ and $A_{ij}$ is the $(ij)^{th}$ element of matrix $A$. The dot product ${\bf{a_{i}}}\cdot{\bf{a_{j}^{T}}}$ is equal to the number of common neighbors of vertices $i$ and $j$, and $A_{ij}=1$ if and only if there is an edge between $i$ and $j$ in $E$.
 For example, if $i$ and $j$ are connected and have 3 common neighbors, then $P_{ij}=5$. Since diameter of the graphs considered in this article is less than three, $P_{ij}\geq 1$ for all pairs of vertices $i$, $j$. 

We denote the maximum and minimum phase values at time $t$ by $\phi_{max}^{t} \triangleq \max\limits_{i}\phi_{i}^{t}$ and $\phi_{min}^{t} \triangleq \min\limits_{i}\phi_{i}^{t}$, where $\phi_{i}^{t}$ is a phase of oscillator $i$ at time $t$.
Let $D_{t}$ be defined as a maximum phase difference between two oscillators at time $t$ $(t\geq 0)$, i.e. 
\begin{equation*}
D_{t}~\triangleq~ \phi_{max}^{t}-\phi_{min}^{t},
\end{equation*}
then $\phi_{min}^{t} \leq \phi_{i}^{t} \leq \phi_{max}^{t}$ ($\forall i = 1, \dots, n$). In other words, each phase lies between the minimal and maximal phases $\phi_{min}^{t}$ and $\phi_{max}^{t}$.
Maximum initial (at time $t=0$) pairwise phase difference is denoted by $D_{0}$:  
\begin{equation*}
D_{0}~=~ \phi_{max}^{0}-\phi_{min}^{0}.
\end{equation*}
If we can show that the maximum phase difference is always bounded, i.e. if $D_{t} \leq D \; \forall \;t \geq 0$, where $D$ is a constant satisfying $D_{0} \leq D <\infty$, then the trajectories will be also bounded since the phase average remains the same (for system \eqref{eq:sys2}: $\dot{\phi_{1}}+\dots+\dot{\phi_{n}}=0$). The PIS, therefore, is defined through the maximum phase difference that is bounded by the value of $D$:
\begin{equation}
\text{PIS}\triangleq\{\vec{\phi}\in \mathbb{R}^{n}: \; \max_{i,j}|\phi_i-\phi_j|\leq D, \; \sum\limits_{i=1}^{n}\phi_{i}=0\},
\end{equation}
which is obviously a compact.

We now formulate a general sufficient condition that guarantees that the maximum phase difference is always bounded by a constant $D$ and thus the trajectories are also bounded.

{\bf{Proposition 1}} {\it{ If $D$ is a constant satisfying $D_{0} \leq D <\infty$, and for all times $t\geq 0$ such that $D_{t}=\phi_{max}^{t}-\phi_{min}^{t} = D$, the following condition is satisfied:
\begin{equation} 
\label{eq:prop1}
\begin{split}
&\dot{\phi}_{k}^{t}-\dot{\phi}_{l}^{t} = \bar{\omega}_{k}-\bar{\omega}_{l} \\ 
&- \sum\limits_{i\in N_{k}}\frac{K_{ik}}{n}\cdot \sin(\phi_{k}^{t} - \phi_{i}^t) -\sum\limits_{j\in N_{l}}\frac{K_{jl}}{n}\cdot \sin(\phi_{j}^{t} - \phi_{l}^t) \leq 0,  \\ 
\end{split} 
\end{equation}
for every two oscillators $k$ and $l$ such that $\phi_{k}^{t}=\phi_{max}^{t}$ and $\phi_{l}^{t}=\phi_{min}^{t}$, then the maximum phase difference is bounded by $D$, i.e. $D_{t}\leq D$ for all $t\geq 0$, trajectories of system \eqref{eq:sys2} are bounded, and system \eqref{eq:sys2} achieves frequency synchronization.}}

\begin{proof}
Condition \eqref{eq:prop1} says that when the maximum phase difference achieves value $D$, it can not grow anymore and thus does not exceed $D$. This implies that the trajectories of system \eqref{eq:sys2} are bounded in $\mathbb R^{n}$ since the phase average is always equal to zero. Further, function $V(\vec{\phi})$ is well-defined in $\mathbb R^{n}$ and we can apply LaSalle's Invariance Theorem to guarantee that each solution of \eqref{eq:sys2} approaches the nonempty set $\{\dot{V} \equiv 0\}=\{\dot{\phi}_{i}=0, 1\leq i\leq n\}$, and system \eqref{eq:sys2} achieves frequency synchronization.
\end{proof}
It is possible that when $\phi_{max}^{t}-\phi_{min}^{t} = D$, several oscillators have phase values equal to $\phi_{max}^{t}$ or $\phi_{min}^{t}$. In this case condition \eqref{eq:prop1} should be satisfied for each pair of oscillators with a phase difference equal to $D$. 

Condition \eqref{eq:prop1} is very general by itself and cannot be directly applied to ensure boundedness of the trajectories and frequency synchronization of a given system.
In the next two subsections we derive two conditions that can be easily verified for each given system and guarantee that condition \eqref{eq:prop1} of Proposition 1 is satisfied. In particular, in Subsection B we derive an analytic condition for the case of equal coupling strengths, and in Subsection C we formulate an optimization problem for the case of non-equal coupling strengths.


An alternative line of works \cite{c3}-\cite{c25} focuses on results that guarantee existence of a locally stable equilibrium manifold for system \eqref{eq:sys2}. These local results, however, cannot guarantee synchronization for any given values of the initial phases (different from the equilibrium phases).
We finish this subsection with an example that demonstrates that existence of a locally stable equilibrium for system \eqref{eq:sys2} does not imply synchronization of this system for all possible initial phases. Therefore, existence of an equilibrium is not a sufficient condition of synchronization for all initial phases.

\begin{figure}[t]
      \centering
      \includegraphics[scale=0.5]{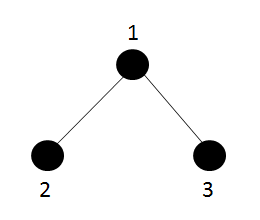}
      \caption{ Topology in Example 1: a star graph with 3 nodes}
      \label{fig2}
   \end{figure}

\begin{figure}[t!]
      \centering
      \includegraphics[scale=0.35]{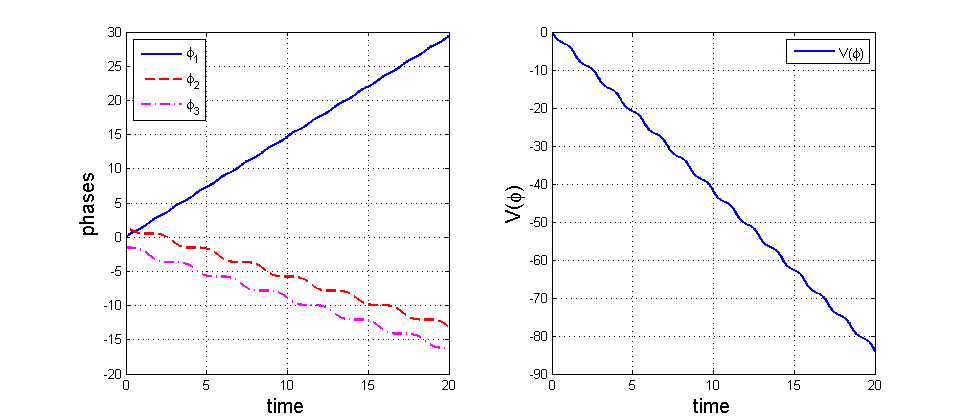}
      \caption{Trajectories $\vec{\phi}(t)$, $t\geq 0$ (left), and $V(\vec{\phi})$ (right) in Example 1}
      \label{fig3}
   \end{figure}

{\bf{Example 1}} In this example three oscillators are connected as shown on Fig. \ref{fig2}, i.e. they form a star graph with three nodes.
We assume that $\bar{\omega}_{1}=2-\epsilon$, $\bar{\omega}_{2}= \bar{\omega}_{3}=-1+~\frac{\epsilon}{2}$, where $\epsilon$ is a small positive constant, and all coupling strengths are equal: $K_{12}=K_{13}=3$.
It is easy to verify that this system possesses a locally stable equilibrium: $\phi_{1} = \frac{2}{3}\sin^{-1}(1-\frac{\epsilon}{2})$, $\phi_{2} = \phi_{3} = -\frac{1}{3}\sin^{-1}(1-\frac{\epsilon}{2}).$ However, there are initial phases $\phi_{1}^{0}, \phi_{2}^{0}$ and $\phi_{3}^{0}$ for which the system does not achieve synchronization. On Fig. \ref{fig3} the behavior of oscillators is demonstrated for $\phi_{1}^{0}= 0$, $\phi_{2}^{0}=\pi/2$ and $\phi_{3}^{0}=-\pi/2$, and for $\epsilon=0.1$. A graph on the right side of Fig. \ref{fig3} is a graph of the Lyapunov function $V(\vec{\phi})$. This function decreases but is not bounded in this example.

\subsection{Analytic Synchronization Condition for System \eqref{eq:sys2} with Equal Coupling Strengths}
In this subsection we consider a special case of system \eqref{eq:sys2} when the coupling strengths are equal for all connected oscillators, i.e. we study the following system:
\begin{equation}
\label{eq:B}
\dot \phi_{i}=\bar{\omega}_{i}+\frac{K}{n}\sum_{j \in N_{i}} \sin(\phi_{j}-\phi_{i}).
\end{equation}

The main result of this subsection is Theorem 1 which contains requirements on the initial phases and the coupling strength such that condition \eqref{eq:prop1} of Proposition 1 is satisfied and therefore system \eqref{eq:B} achieves frequency synchronization. 

$\bf{Theorem \;1}$ {\it {If $D$ is a constant satisfying $0<D_{0}\leq D <\pi$, and
\begin{equation}
\label{eq:thm1}
K \geq \frac{n\cdot |\bar{\omega}_{i} - \bar{\omega}_{j}|}{P_{ij}\cdot \sin D} 
\end{equation}
for all $i,j = 1, \dots, n$, then  $D_{t} \leq D \; \forall \;t \geq 0$ for the system \eqref{eq:B} in which the underlying topology is a graph with diameter~$\leq 2$, and this system achieves frequency synchronization.}} 

\begin{proof}
Assume that at time moment $T \geq0$, the value of $D_{T}$ is equal to $D$ and before this moment it never exceeded $D$, i.e. $D_{t}\leq D \; \forall t\leq T$. We will show that under the conditions of this theorem, the maximum phase difference does not start to increase at time $T$ by showing that requirement \eqref{eq:prop1} of Proposition 1 is satisfied. This will guarantee that the maximum phase difference $D_{t}$ will be always bounded by $D$.

Condition \eqref{eq:prop1} must be satisfied for every two oscillators $k$ and $l$ with $\phi_{k}^{T}=\phi_{max}^{T}$ and $\phi_{l}^{T}=\phi_{min}^{T}$: 
\begin{equation*}
\begin{split}
&\dot{\phi}_{k}^{T}-\dot{\phi}_{l}^{T} = \bar{\omega}_{k}-\bar{\omega}_{l} \\
&- \frac{K}{n}\sum\limits_{i\in N_{k}}\sin(\phi_{k}^{T} - \phi_{i}^T)-\frac{K}{n}\sum\limits_{j\in N_{l}}\sin(\phi_{j}^{T} - \phi_{l}^T) \leq 0. \\
\end{split}
\end{equation*} 
This condition will be satisfied if $K \geq \frac{n\cdot |\bar{\omega}_{k} - \bar{\omega}_{l}|}{P_{kl}\cdot \sin D}$ and if we can show that:
\begin{equation}
\label{eq:th1_proof}
\sum\limits_{i\in N_{k}}\sin(\phi_{k}^{T} - \phi_{i}^T)+\sum\limits_{j\in N_{l}}\sin(\phi_{j}^{T} - \phi_{l}^T) \geq P_{kl}\cdot \sin D.
\end{equation}
Because $\phi_{k}^{T}$ and $\phi_{l}^{T}$ are respectively the maximum and minimum phase values at time $T$ (see Fig. \ref{fig4}):
\begin{equation*}
\begin{split}
&0 \leq \phi_{k}^{T}-\phi_{i}^{T}\leq D<\pi,\\
&0 \leq \phi_{j}^{T}-\phi_{l}^{T}\leq D<\pi,
\end{split}
\end{equation*}
where $1\leq i,j\leq n$. Therefore, each summand in the left side of the inequality \eqref{eq:th1_proof} is nonnegative.

If vertices $k$ and $l$ are connected by an edge, both sums contain $\sin(\phi_{k}^{T}-\phi_{l}^{T})=\sin D$, and thus the left side of \eqref{eq:th1_proof} contains $2\sin D$.

Assume now that vertices $k$ and $l$ have a common neighbor ~-~ vertex $m$.
Then, the left side of inequality \eqref{eq:th1_proof} contains the following sum:
\begin{equation*}
\begin{split}
&\sin(\phi_{k}^{T}-\phi_{m}^{T})+\sin(\phi_{m}^{T}-\phi_{l}^{T})\\
&=2\sin\Bigl(\frac{\phi_{k}^{T}-\phi_{l}^{T}}{2}\Bigl)\cdot \cos\Bigl(\frac{\phi_{k}^{T}+\phi_{l}^{T}-2\phi_{m}^{T}}{2}\Bigl)\\
&\geq2\sin\Bigl(\frac{\phi_{k}^{T}-\phi_{l}^{T}}{2}\Bigl)\cdot \cos\Bigl(\frac{\phi_{k}^{T}-\phi_{l}^{T}}{2}\Bigl)\\
&=\sin(\phi_{k}^{T}-\phi_{l}^{T})=\sin D.
\end{split}
\end{equation*} 
Inequality above holds because $\sin\bigl(\frac{\phi_{k}^{T}-\phi_{l}^{T}}{2}\bigl)>0$, and 
\begin{equation*}
-\frac{\pi}{2}<-\frac{D}{2}\leq \frac{\phi_{k}^{T}+\phi_{l}^{T}-2\phi_{m}^{T}}{2}\leq \frac{D}{2}<\frac{\pi}{2},
\end{equation*}
so that 
\begin{equation*}
 \cos\Bigl(\frac{\phi_{k}^{T}+\phi_{l}^{T}-2\phi_{m}^{T}}{2}\Bigl)\geq \cos\Bigl(\frac{\phi_{k}^{T}-\phi_{l}^{T}}{2}\Bigl)=\cos\Bigl(\frac{D}{2}\Bigl).
\end{equation*}

   \begin{figure}[t]
      \centering
      \includegraphics[scale=0.5]{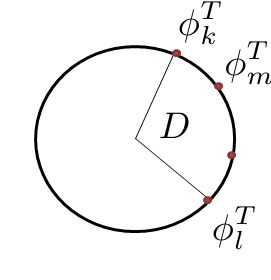}
      \caption{ Oscillators $\phi_{k}^{T}$ and $\phi_{l}^{T}$ with the maximum and minimum phases, respectively}
      \label{fig4}
   \end{figure}

Therefore, the left side of \eqref{eq:th1_proof} contains a sum that is greater or equal than $\sin D$ for each common neighbor $m$ of vertices $k$ and $l$. In addition, if $k$ and $l$ are connected by the edge, there is a term $2\sin D$ in the left side of \eqref{eq:th1_proof}, and thus, inequality \eqref{eq:th1_proof} holds. This proves that condition \eqref{eq:prop1} of Proposition 1 is satisfied under the theorem's conditions.
\end{proof}

{\bf{Remark 1}} If $D_{0}\leq \frac{\pi}{2}$, the smallest value of bound \eqref{eq:thm1} will be achieved for $D=\frac{\pi}{2}$. When $\frac{\pi}{2}<D_{0}<\pi$, bound \eqref{eq:thm1} takes its smallest value if $D=D_{0}$. 

{\bf{Remark 2}} In the case of a complete graph, $P_{ij}=n$ for each pair $i$, $j$ of vertices, and the sufficient condition on $K$ is the following: $K \geq \frac{|\bar{\omega}_{i} - \bar{\omega}_{j}|}{\sin D_{0}}$ for all $ i, j$. This bound coincides with the bound obtained in \cite{c3} for a complete graph.

{\bf{Remark 3}} When diameter of a graph is larger than two, Theorem 1 cannot be applied in general, and condition \eqref{eq:prop1} of the Proposition 1 can be violated. For instance, if the distance between vertices $k$ and $j$ is more than two, then in condition \eqref{eq:prop1} both sums may be equal to zero: $\sin(\phi_{k}^{t}-\phi_{i}^{t})=0 \; \forall i~\in~N_{k}$,   $\sin(\phi_{j}^{t}-\phi_{l}^{t})=0 \; \forall j \in N_{l}$, and $(\dot{\phi}_{k}^{t}-\dot{\phi}_{l}^{t})>0$ if $\bar{\omega}_{k}>\bar{\omega}_{l}$.
However, Theorem 1 can be applied to the graphs with a diameter more than two if every two oscillators with a shortest path between them of a length more than two, have equal frequencies. In this case condition \eqref{eq:prop1} is always satisfied for such two oscillators. Indeed, if $\bar{\omega}_{k}=\bar{\omega}_{l}$ and $(\phi_{k}^{t}-\phi_{l}^{t})=(\phi_{max}^{t}-\phi_{min}^{t})=D<\pi$, then $\dot{\phi}_{k}^{t}-\dot{\phi}_{l}^{t}\leq 0$, because $\sin(\phi_{k}^{t} - \phi_{i}^t)\geq 0$ and $\sin(\phi_{j}^{t} - \phi_{l}^t)\geq 0$ for all $i \in N_{k}$ and $j \in N_{l}$.

\subsection{Optimization Approach for System \eqref{eq:sys2} with Non-equal Coupling Strengths}
In this subsection the equal coupling strength assumption is relaxed.
 Instead of one coupling parameter $K$ as was in the previous subsection, there are now $|E|$ coupling parameters $K_{ij}$, where $|E|$ is the cardinality of the graph's edge set $E$. Similarly to condition \eqref{eq:thm1} in the Theorem 1, we will find bounds on the coupling strengths $K_{ij}$ to guarantee frequency synchronization of system \eqref{eq:sys2}, but instead of providing an analytic condition \eqref{eq:thm1}, we will formulate an optimization problem whose solution contains the coupling strengths $K_{ij}$ that guarantee \eqref{eq:prop1} and are sufficient for synchronization.

While in the Theorem 1 the goal was to find the minimum value of the coupling parameter $K$ that guarantees synchronization, minimizing the sum of all coupling strengths $\sum\limits_{ij\in E}K_{ij}$ will be the goal for the case of non-equal coupling strengths\footnote{Objective function used in this article, therefore, is the $l^{1}$-norm of a vector of all coupling strengths. Other options could be employed, for example, the maximum norm $l^{\infty}$ could be used that corresponds to minimizing the maximum coupling strength.}.


In condition \eqref{eq:prop1} we assume that $\phi_{k}^{t}=\phi_{max}^{t}$, $\phi_{l}^{t}=\phi_{min}^{t}$ and $\phi_{k}^{t}-\phi_{l}^{t}=D$. Since $D<\pi$, all values of $\sin()$ functions in each sum of \eqref{eq:prop1} are nonnegative. Instead of condition \eqref{eq:prop1} we will consider a more strict condition on the coupling parameters, where we keep only summands corresponding to the neighbor oscillators of both oscillators $k$ and $l$:
\begin{equation}
\label{eq:ineq}
\begin{split}
&\dot{\phi}_{k}^{t}-\dot{\phi}_{l}^{t} \leq \bar{\omega}_{k}-\bar{\omega}_{l}-\frac{2K_{kl}}{n}\cdot \sin D\\
&-\sum\limits_{m\in N_{kl}}\Bigl(\frac{K_{km}}{n}\cdot\sin(\phi_{k}^{t} - \phi_{m}^t)+\frac{K_{lm}}{n}\cdot\sin(\phi_{m}^{t} - \phi_{l}^t)\Bigl) \leq 0,  \\ 
\end{split} 
\end{equation}
where $N_{kl}=N_{k} \cap N_{l}$ -- is the set of common neighbors of oscillators $k$ and $l$. If there is no edge $kl$ between oscillators $k$ and $l$, then $K_{kl}=0$ in \eqref{eq:ineq}.
We will introduce constraints that do not contain phases and guarantee that condition \eqref{eq:ineq} (and \eqref{eq:prop1} as well) is satisfied for all phase values. Optimization problem, whose $|E|$ variables are the coupling strengths $K_{ij}$ ($ij \in E$), that allowed to take nonnegative values, is formulated as follows:
\vspace{10pt}

\hspace{-8pt}{\bf{minimize:}} \hspace{10pt} $\sum\limits_{ij\in E}K_{ij},$ 
\begin{equation}
\label{eq:constraints}
\begin{split} 
&\text{{\bf{\hspace{-10pt}subject to: \hspace{5pt}    }}}|\bar{\omega}_{k}-\bar{\omega}_{l}| -\frac{2K_{kl}}{n}\cdot \sin D\\
&-\frac{\sin D}{n}\sum\limits_{m \in N_{kl}}\Bigl(\delta_{m}\cdot K_{km}+(1-\delta_{m})\cdot K_{lm}\Bigl)\leq 0,
\end{split}
\end{equation}
where $1\leq k,l\leq n$, and each $\delta_{m}$ may take values $\{0,1\}$. Since either $\delta_{m}$ or $(1-\delta_{m})$ takes a zero value, variables $K_{km}$ and $K_{lm}$ do not appear together in each constraint.

For each possible combination of values of $\delta_{m}$ there is a corresponding constraint, and, therefore for each pair of oscillators $k$ and $l$ there are $2^{|N_{kl}|}$ constraints in the optimization problem, where $|N_{kl}|$ is the number of common neighbors of oscillators $k$ and $l$. For example, suppose that oscillators $k$ and $l$ are connected and have a single common neighbor $m$, then optimization problem \eqref{eq:constraints} will contain two constraints for oscillators $k$ and $l$:
\begin{equation*}
\begin{split}
&|\bar{\omega}_{k}-\bar{\omega}_{l}| -\frac{2K_{kl}}{n}\cdot \sin D-\frac{\sin D}{n} K_{km}\leq 0 \text{, when }\delta_{m}=1,\\
&|\bar{\omega}_{k}-\bar{\omega}_{l}| -\frac{2K_{kl}}{n}\cdot \sin D-\frac{\sin D}{n} K_{lm}\leq 0 \text{, when }\delta_{m}=0.
\end{split}
\end{equation*}
If, for example, oscillators $k$ and $l$ are not connected and have two common neighbors $m_{1}$, $m_{2}$, then there will be four constraints for $k$ and $l$:
$$|\bar{\omega}_{k}-\bar{\omega}_{l}| -\frac{\sin D}{n}\Bigl(K_{km_{1}}+K_{km_{2}}\Bigl)\leq 0, \text{ } (\delta_{m_{1}}=\delta_{m_{2}}=1),$$
$$|\bar{\omega}_{k}-\bar{\omega}_{l}| -\frac{\sin D}{n}\Bigl(K_{lm_{1}}+K_{km_{2}}\Bigl)\leq 0, \text{ } (\delta_{m_{1}}=0, \delta_{m_{2}}=1),$$
$$|\bar{\omega}_{k}-\bar{\omega}_{l}| -\frac{\sin D}{n}\Bigl(K_{km_{1}}+K_{lm_{2}}\Bigl)\leq 0, \text{ } (\delta_{m_{1}}=1, \delta_{m_{2}}=0),$$
$$|\bar{\omega}_{k}-\bar{\omega}_{l}| -\frac{\sin D}{n}\Bigl(K_{lm_{1}}+K_{lm_{2}}\Bigl)\leq 0, \text{ } (\delta_{m_{1}}=\delta_{m_{2}}=0).$$
Thus, optimization problem \eqref{eq:constraints} contains in total  $\sum\limits_{1\leq k<l\leq n}2^{|N_{kl}|}$ constraints. 
Although the number of constraints can be exponential in number of oscillators $n$, for some types of graphs it is polynomial in $n$. For example, for the graphs with star-tree topology, each pair of oscillators has at most one common neighbor, and thus, not more than two corresponding constraints. 

{\bf{Remark 4}} If all coupling strengths are required to be equal in optimization problem \eqref{eq:constraints}, then its solution is bound \eqref{eq:thm1} from the Theorem 1. Indeed, when all coupling strengths are equal, then $K_{kl}=K_{km}=K_{lm}$ in the constraint of \eqref{eq:constraints} for $\bar{w}_{k}$ and $\bar{w}_{l}$, $\sum\limits_{m \in N_{kl}}\Bigl(\delta_{m}\cdot K_{km}+(1-\delta_{m})\cdot K_{lm}\Bigl)=|N_{kl}|$, and the constraint becomes: $K\geq\frac{n\cdot |\bar{\omega}_{k}-\bar{\omega}_{l}|}{P_{kl}\cdot \sin D}$.

\begin{table*}[t]

\caption{Synchronization conditions in our comparative analysis}
\begin{center}
    \begin{tabular}{ | c | c | c |}
    \hline

    & \bf{Bound on Coupling Strength} & \bf{Constraint on Initial Phases} \\ \hline
&&\\
  Analytic condition (Theorem 1) & $K \geq \frac{n\cdot |\bar{\omega}_{i} - \bar{\omega}_{j}|}{P_{ij}\cdot \sin D} \; (\forall i,j)$ & $ D_{0}<\pi$  \\
 && \\
 \hline
&&\\
  
   Numerical condition (Theorem 2)&  Solution to \eqref{eq:constraints}  & $D_{0}<\pi$\\

&&\\
    \hline
 && \\
    Condition from \cite{c3}& $K> \frac{2n\cdot\left\|B_{c}^{T}\bar{\omega}\right\|_{2}}{\lambda_{2}\cdot \pi \cdot \text{sinc}(\gamma_{max})  }$ & $D_{0}<\pi; \; \left\|B_{c}^{T}\phi(0)\right\|_{2} < \pi$  \\ 
 && \\
\hline
&&\\
  
   Condition from \cite{c5} &  $K > \frac{\sqrt{2}\sigma(\bar{\omega})}{L_{*}\cdot\sin D}$  & $D_{0}<\pi; \; \sum\limits_{i=1}^{n}\phi_{i}^{0}=0; \; \mathcal{E}_{0}<\frac{D^{2}}{2}<\frac{\pi^{2}}{2}$\\

&&\\
    \hline

&&\\
 &$ K \geq \frac{\sigma(\bar{\omega})\cdot D}{\sqrt{\mathcal{E}_{0}}\cdot L\cdot \sin D}$&\\
    Numerical condition from \cite{c28}  && $\mathcal{E}_{0}< D^{2}<\pi^{2}$ \\
&$K \geq \frac{n\cdot|\bar{\omega}_{k}-\bar{\omega}_{l}|}{\sum\limits_{i \in N_{k}}\sin(\phi_{k}-\phi_{i})+\sum \limits_{j \in N_{l}}\sin(\phi_{j}-\phi_{l})}$& \\
 \hline
   
     \end{tabular}

\label{table_example}
\end{center}
\end{table*}  

We will now show that solution to this optimization problem satisfies conditions \eqref{eq:ineq} for all possible phase values.

{\bf{Theorem 2}} {\it{Solution to the optimization problem \eqref{eq:constraints} satisfies conditions \eqref{eq:prop1}, and system \eqref{eq:sys2} achieves frequency synchronization.}}

\begin{proof}
Suppose that $K_{ij}^{*}$, where $ij \in E$  is a solution of the optimization problem \eqref{eq:constraints}. We are going to  show that condition \eqref{eq:ineq} is satisfied for two arbitrary oscillators $k$ and $l$ with $\phi_{k}^{t}-\phi_{l}^{t}=D$. This would imply that condition \eqref{eq:prop1} is also satisfied since condition \eqref{eq:ineq} is more restrictive than \eqref{eq:prop1}.

 For arbitrary phases $\phi_{m}^{t}$, such that $\phi_{l}^{t}\leq \phi_{m}^{t}\leq \phi_{k}^{t}$ for all $m \in N_{kl}$, from condition \eqref{eq:ineq}: 
\begin{equation*}
\begin{split}
&\bar{\omega}_{k}-\bar{\omega}_{l}-\frac{2K_{kl}^{*}}{n}\cdot \sin D\\
&-\sum\limits_{m\in N_{kl}}\Bigl(\frac{K_{km}^{*}}{n}\cdot\sin(\phi_{k}^{t} - \phi_{m}^t)+\frac{K_{lm}^{*}}{n}\cdot\sin(\phi_{m}^{t} - \phi_{l}^t)\Bigl)\\
& \leq\bar{\omega}_{k}-\bar{\omega}_{l}-\frac{2K_{kl}^{*}}{n}\cdot \sin D\\
&-\sum\limits_{m \in N_{kl}}\frac{\min(K_{km}^{*},K_{lm}^{*})}{n}\cdot \bigl(\sin(\phi_{k}^{t} - \phi_{m}^t)+\sin(\phi_{m}^{t} - \phi_{l}^t)\bigl)\\
&\leq \bar{\omega}_{k}-\bar{\omega}_{l}-\frac{2K_{kl}^{*}}{n}\cdot \sin D-\frac{\sin D}{n}\sum\limits_{m \in N_{kl}}\min(K_{km}^{*},K_{lm}^{*}),
\end{split}
\end{equation*}
because $(\phi_{k}^{t}-\phi_{m}^{t}) \in [0,D]$, $(\phi_{m}^{l}-\phi_{l}^{t}) \in [0,D]$, and $(\phi_{k}^{t}-\phi_{m}^{t})+(\phi_{m}^{t}-\phi_{l}^{t})=D$.

Now we can observe that for the right side of the last inequality there exists a constraint in \eqref{eq:constraints} that guarantees that the right side is non-positive.
If, for example, $\min(K_{km}^{*},K_{lm}^{*})=K_{km}^{*}$, then corresponding constraint in \eqref{eq:constraints} has $\delta_{m}=1$, otherwise $\delta_{m}=0$.
\end{proof}

We finish this section with an example for which we found values of the coupling strengths that are sufficient for synchronization: first, under the condition that all coupling strengths must be equal and using the Theorem 1, and then, assuming that the strengths are allowed to be non-equal and solving  the optimization problem \eqref{eq:constraints}.

 \begin{figure}[b]
      \centering
      \includegraphics[scale=0.5]{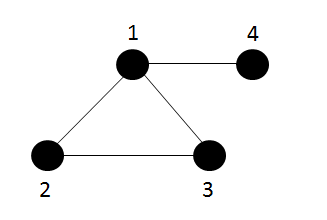}
      \caption{Connections between oscillators in Example 2}
      \label{fig5}
   \end{figure}

{\bf{Example 2}} In this example we consider four oscillators connected as shown on Fig. \ref{fig5} and with following frequencies: $\bar{\omega}_{1}=-0.1$, $\bar{\omega}_{2}=0$, $\bar{\omega}_{3}=0.3$, $\bar{\omega}_{4}=-0.2$. There are four edges in this graph, i.e. four coupling strengths $K_{ij}$, and thus four variables in problem \eqref{eq:constraints}. Notice, that $P_{12}=P_{13}=P_{23}=3$, $P_{14}=2$, and $P_{24}=P_{34}=1$. If we assume that all the coupling strengths are equal, then by the Theorem 1 from previous subsection, sufficient for synchronization value of the coupling strength is: $K=0.5\cdot n = 2$ (from the inequality for pair $34$). Then, the sum of all coupling strengths is $4\cdot2=8$.

If we let the coupling strengths be different for the different edges, the optimization problem has a solution: $K_{12}=0.8$, $K_{13}=2$, $K_{23}=0.2$, and $K_{14}=2$. Now the sum of the coupling strengths is $5$. Optimization problem for this example contains eleven inequality constraints (besides the constraints $K_{ij} \geq 0$).


For optimization we used Matlab's R2012a $\it{fmincon}$ function with default options.



\section{Numerical Simulations}

In this section we present the results of simulations performed to demonstrate that for the graphs of diameter two, synchronization condition formulated in Theorem 1 is a less restrictive condition compared to existing ones. Since Theorem 1 guarantees existence of a Positively Invariant Set and frequency synchronization of system \eqref{eq:sys2}, we compared our bound with the similar conditions that also guarantee existence of a PIS and frequency synchronization. To the best of our knowledge, there are three such conditions: Theorem 4.6 from \cite{c3},  results from \cite{c5}, and conditions (analytic and numerical) in \cite{c28}. Therefore, we did not include into comparison analysis conditions from \cite{c3}, \cite{c24} and \cite{c25} that only provide existence of an equilibrium and local stability.
The numerical condition of \cite{c28} is less restrictive then the analytic condition of the same article, and we here consider only the former one.
In addition, we added to our comparison analysis a numerical synchronization condition from Theorem 2, which allows the coupling strengths to be different, and for each given example we calculated an average coupling strength of the solution to \eqref{eq:constraints}.

\begin{figure*}[t]
\begin{center}
  \begin{subfigure}[t]{0.4\textwidth}
      \includegraphics[width=\textwidth]{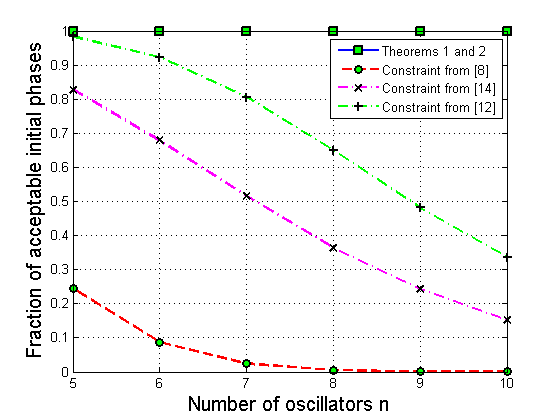}
      \caption{Fractions of random samples of initial phases that satisfy initial phase constraints }  
      \label{exp3_a}    
\end{subfigure}
\begin{subfigure}[t]{0.4\textwidth}
      \includegraphics[width=\textwidth]{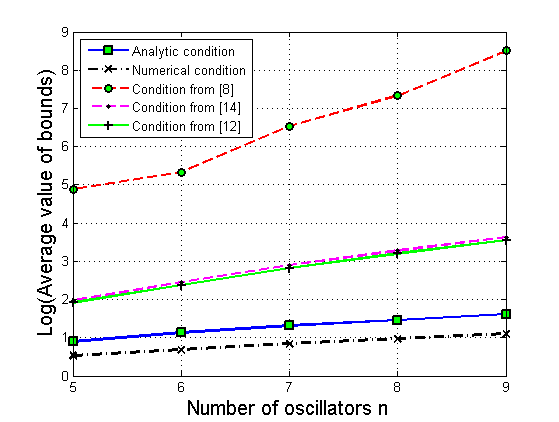}
      \caption{Average values of bounds in logarithmic scale}
      \label{exp3_b}          
\end{subfigure}
\end{center}
\caption{Simulation results of Experiments 1 and 2}
   \end{figure*}

Each of the five synchronization conditions compared in this section consists of a bound on the coupling strength and constraints on the initial phases of oscillators.
 In particular, all synchronization conditions require that the difference between any two initial phases is less than $\pi$ (i.e. $D_{0}<\pi$).
Additionally, synchronization conditions from \cite{c3}, \cite{c28} and \cite{c5} have their own special constraints on the initial phases.
The bounds on the coupling strength and corresponding requirements on the initial phases are summarized in Table~1.

In the simulations we assigned a value of $\max\{\frac{\pi}{2},D_{0}\}$ to the constant $D$ for our synchronization condition, because in this case bound \eqref{eq:thm1} is the least restrictive as mentioned in Remark 1.
In the bound from \cite{c3},  $\lambda_{2}$ is the algebraic connectivity of a given graph, $B_{c} \in \mathbb{R}^{n\times n(n-1)/2}$ is the incidence matrix of the complete graph with $n$ nodes, $\bar{\omega}$ is a vector of frequencies, $\phi(0)$ ~--~ vector of initial phases and  $\gamma_{max}=\max\{\frac{\pi}{2},\left\|B_{c}^{T}\phi(0)\right\|_{2}\}$.
 In the condition from \cite{c5}, $\mathcal{E}_{0}$ is the squared Euclidean norm of a vector of the initial phases:
\begin{equation*}
\mathcal{E}_{0} \triangleq \sum\limits_{i=1}^{n}\bigl(\phi_{i}^{0}\bigl)^{2},
\end{equation*}
$\sigma(\bar{\omega})$ denoted he Euclidean norm of a vector of the intrinsic frequencies deviations:
\begin{equation*}
\sigma(\bar{\omega})\triangleq \sqrt{\sum\limits_{i=1}^{n}(\bar{w}_{i})^{2}},
\end{equation*}
$D$ is a constant whose value is defined as $\max\{\frac{\pi}{2}, \sqrt{2\mathcal{E}_{0}}\}$, and $L_{*}$ is defined as $L_{*}\triangleq \frac{1}{1+diam(G)\cdot |E^{c}(G)|}$, where $diam(G)=2$ is the diameter of a graph $G$ and $|E^{c}(G)|$ is the cardinality of the set $E^{c}(G)$ defined as
\begin{equation*}
E^{c}(G)\triangleq E_{comp} \setminus E(G),
\end{equation*}
where $E_{comp}$ is a set of $\frac{n(n-1)}{2}$ edges of a complete graph with $n$ nodes. 

In the condition from \cite{c28}, $\sigma(\bar{\omega})$ and $\mathcal{E}_{0}$ are defined as in \cite{c5}, $L=L_{*}$, and $D=D_{0}$.

In our analysis we compared the requirements on both, the initial phases, and on the coupling strength of the five synchronization conditions.

{\bf{Experiment 1 (comparison of the constraints on initial phases).}} 
In the first experiment we checked the restrictiveness of constraints on the initial phases of each of five synchronization conditions under consideration. We created $10^{5}$ samples of the initial phases such that each phase was chosen from the $(0,\pi)$ interval. Then, for each sample we subtracted its mean phase value from each phase belonging to this sample. Therefore, the sum of the initial phases was equal to zero, and the maximum phase difference was less than $\pi$ in each sample. We shifted the phase values of each sample by the sample's mean because condition from \cite{c5} requires that $\sum\limits_{i=1}^{n}\phi_{i}^{0}=0$, and other synchronization conditions only depend on the relative values of the initial phases and thus are rotationally invariant.

Next, for each sample we checked if it satisfies the constraints on the initial phases of the synchronization conditions, and for each condition we calculated fractions of samples that satisfy its initial phase requirements.


We repeated this experiment for different numbers of oscillators $n$ in the system: $n=5,\dots,10$ and the experiment's results are shown on Fig. \ref{exp3_a}.
Since our synchronization conditions in Theorems 1 and 2 do not contain any additional requirements on the initial phases, they can be applied for each generated sample of phases, and thus the fraction of acceptable initial phases is equal to one for all values of $n$.

Fractions of acceptable initial phases for conditions from \cite{c3},  \cite{c5} and \cite{c28} monotonically decrease with the number of oscillators $n$ as can be observed on Fig. \ref{exp3_a}. 




{\bf{Experiment 2 (comparison of the bounds on coupling strength).}}
In the second experiment we compared the bounds on the coupling strength. For each fixed number of oscillators $n=5,\dots,9$ we randomly created 1000 graphs with $n$ vertices and of diameter two.  The initial edge set of each graph was empty, and we successively added random edges to it until the diameter was equal to two.
 For each graph we then created a random sample of initial phases, a random sample of frequencies, and calculated the bounds on $K$ for each condition. For the numerical condition in the Theorem 2 we calculated an average value of $K$ for each example. The average values of bounds for each of five conditions under comparison are plotted on Fig. \ref{exp3_b} in logarithmic scale.
In this experiment we sampled values of the frequencies from $(0,1)$ interval, but the relative performance of the bounds does not noticeably change with the interval.

The simulation results of Experiments 1 and 2 show that  for graphs of diameter two our synchronization condition formulated in Theorem 1 is less restrictive in terms of both, initial phases and coupling strength compared to the existing conditions. Additionally, optimization-based condition in Theorem 2 provides a further improvement if the value of its bound is defined as the average coupling strength for each example.

\section{Conclusion}

In this article we employed the notion of a Positively Invariant Set to find a sufficient condition for frequency synchronization of heterogeneous Kuramoto oscillators connected by a graph of diameter two. We showed that an existence of a PIS ensures the boundedness of the trajectories of oscillators, which in turn, provides synchronization.
For the case when the coupling strength is the same for every two connected oscillators, we provided an analytic synchronization condition, and  demonstrated with simulations that this condition is significantly less restrictive than existing ones. For the case when the coupling is allowed to take distinct values for different pairs of oscillators, we formulated an optimization problem whose solution~--~a set of coupling strengths~--~guarantees frequency synchronization.

\end{document}